\documentclass[11pt,a4paper]{article}

\usepackage{amsmath,amssymb,amsthm,mathtools}
\usepackage[margin=1in]{geometry}
\usepackage[colorlinks=true,linkcolor=blue,urlcolor=blue]{hyperref}

\usepackage{graphicx, tikz}
\usetikzlibrary{arrows.meta}

\newtheorem{theorem}{Theorem}[section]
\newtheorem{lemma}[theorem]{Lemma}
\newtheorem{remark}[theorem]{Remark}

\numberwithin{equation}{section}
\numberwithin{figure}{section}

\newcommand{\E}{\mathbb{E}}
\newcommand{\R}{\mathbb{R}}
\newcommand{\supp}{\operatorname{supp}}
\newcommand{\TV}{d_{\mathrm{TV}}}
\DeclareMathOperator{\Tr}{Tr}
\DeclareMathOperator{\Var}{Var}

\newcommand{\eps}{\epsilon}

\newcommand{\email}[1]{\href{mailto:#1}{\texttt{#1}}}

\let\oldthebibliography\thebibliography
\let\endoldthebibliography\endthebibliography

\renewenvironment{thebibliography}[1]
  {
    \oldthebibliography{#1}
    \footnotesize
    \setlength{\itemsep}{1pt}
    \setlength{\parskip}{0pt}
    \setlength{\parsep}{0pt}
  }
  {
    \endoldthebibliography
  }

\title{Mesoscopic eigenvalue statistics for random regular graphs: \\ a case study}
\author{Jaehun Lee\textsuperscript{*}}
\date{\today}

\begin{document}

\maketitle

\begingroup
\renewcommand{\thefootnote}{\fnsymbol{footnote}}
\footnotetext[1]{
ISTA, Am Campus 1, 3400 Klosterneuburg, Austria.
Email: \email{jaehun.lee@ist.ac.at}.
}
\endgroup

\begin{abstract}
We study mesoscopic fluctuations of linear eigenvalue statistics for random $d$-regular graphs on $N$ vertices with fixed $d\ge 3$.
Revisiting the approach developed by Dumitriu--Johnson--Pal--Paquette~\cite{MR3078290} and Johnson~\cite{Johnson2015}, we establish a mesoscopic central limit theorem in the bulk and identify its universal limiting variance.
Our result applies to test functions with smooth compactly supported Fourier transforms and to a specific range of spectral scales $(\log N)^{-1}\ll\eta\ll 1$.
In particular, this case study allows the mesoscopic scale to tend to zero arbitrarily slowly.
\end{abstract}

\section{Introduction}
Fix an integer $d\geq 3$.
Let $\mathcal{G}_N$ be chosen uniformly at random from the set of simple $d$-regular graphs on the vertex set $\{1,\ldots,N\}$.
Let $A\equiv A_N$ be the adjacency matrix of $\mathcal{G}_N$.
Equivalently, identifying each graph with its adjacency matrix, one can view $A$ as a random matrix uniformly distributed on the set
\begin{equation*}
	\Big\{A\in\{0,1\}^{N\times N}: A=A^{\mathsf T},\; A_{ii}=0,\;
	\sum_{j}A_{ij}=d \text{ for every }i\in\{1,\ldots,N\}\Big\}.
\end{equation*}
This is the \emph{uniform model} of random regular graphs.
Define the rescaled matrix $H\equiv H_N$ by
\begin{equation*}
	H:=\frac{A}{\sqrt{d-1}}.
\end{equation*}
The eigenvalues of $H$ are denoted by $\{\lambda_{i}\}_{i=1}^{N}$.

For a suitable class of analytic test functions $g$, Johnson~\cite{Johnson2015} showed that the linear eigenvalue statistic
\begin{equation*}
	\Tr g(H) = \sum_{i} g(\lambda_{i})
\end{equation*}
converges in distribution to an infinitely divisible random variable after appropriate deterministic centering.
An analogous result had previously been established for the \emph{permutation model} of random regular graphs by Dumitriu--Johnson--Pal--Paquette~\cite{MR3078290}. The permutation model is constructed by summing a fixed number of independent uniformly random permutation matrices together with their transposes.

This paper concerns a \emph{mesoscopic} version of this eigenvalue fluctuation problem.
For a fixed bulk energy $E\in(-2,2)$ and a spectral scale parameter $\eta\equiv\eta_{N}>0$ such that $\eta\to 0$ and $\eta\gg N^{-1}$, we consider the statistic
\begin{equation}\label{eq: 73}
	\Tr g\left(\frac{H-E}{\eta}\right).
\end{equation}
The aforementioned results in \cite{MR3078290,Johnson2015} correspond to the macroscopic case $\eta\sim 1$.
While these macroscopic fluctuations are described by non-Gaussian infinitely divisible distributions, mesoscopic eigenvalue statistics of various random matrix models exhibit universal Gaussian fluctuations.
In the real symmetry class, the universal limiting variance of \eqref{eq: 73} with a bulk energy takes the form
\begin{equation}\label{eq: 77}
	V_{g} = \frac{1}{2\pi^2} \iint_{\R^2}\left(\frac{g(x)-g(y)}{x-y}\right)^2\mathrm{d}x\mathrm{d}y
	= \frac{1}{2\pi^2} \int_{\R}|t|\,|\widehat g(t)|^2\,dt.
\end{equation}
Such results are known as mesoscopic central limit theorems (CLTs).
The literature on the universality of mesoscopic CLTs for linear eigenvalue statistics is extensive \cite{MR3678478,MR3914908,LS20,MR4255226,LLS24,Vova25a,MR1678012,MR1689027,LS15+,MR3459158,MR3771744,MR3865662,MR4216520,MR4705298}, covering classical ensembles such as Wigner matrices as well as certain classes of random matrices with dependent entries \cite{EL26+,MR4893421,MR4403963}.

For fixed degree $d\geq 3$, the empirical spectral distribution of $H$ converges to the rescaled Kesten--McKay law~\cite{MR109367,McKay1981}. 
At finer scales, local laws and eigenvector delocalization are established throughout the spectrum~\cite{MR3962004,HuangYau2024}. 
GOE universality has subsequently been established for the extreme nontrivial eigenvalues~\cite{HuangMcKenzieYau2024}, and for bulk eigenvalue statistics~\cite{HeHuangWang2026,OpenAI2026RRG}.
These microscopic universality results do not directly imply a mesoscopic CLT, which concerns fluctuations over shrinking spectral windows containing a diverging number of eigenvalues.
In this paper, we establish a mesoscopic CLT for the specific case $(\log N)^{-1}\ll\eta\ll 1$ by revisiting the combinatorial approach developed by Dumitriu--Johnson--Pal--Paquette~\cite{MR3078290} and Johnson~\cite{Johnson2015}. 

Their approach uses a Chebyshev series of the test function $g$ to reduce the problem to counting a certain type of walk.
To adapt this reduction to our mesoscopic setting, we impose a condition on the Fourier transform of $g$.
We use the following Fourier transform convention:
\begin{equation*}
	\widehat{g}(t)=\int_{\R}e^{-\mathrm{i}tx}g(x)\,\mathrm{d}x.
\end{equation*}
Our result on the mesoscopic CLT for fixed-degree random regular graphs is stated as follows.

\begin{theorem}\label{thm:main}
Fix a bulk energy $E\in(-2,2)$. Suppose that the spectral scale $\eta>0$ satisfies
\begin{equation*}
	(\log{N})^{-1} \ll \eta \ll 1,
\end{equation*}
i.e.~$\eta\to 0$ and $\eta\log{N}\to\infty$.
Let $g$ be real-valued and Schwartz, and assume $\widehat{g}\in C_c^\infty(\R)$.
Let $V_{g}$ be as in \eqref{eq: 77}.
Then,
the centered eigenvalue statistic
\begin{equation*}
	\Tr g\left(\frac{H-E}{\eta}\right) - \E\Tr g\left(\frac{H-E}{\eta}\right)
\end{equation*}
converges in distribution to a Gaussian random variable with mean zero and variance $V_{g}$.
\end{theorem}

\begin{remark}[Macroscopic versus mesoscopic fluctuations]
For fixed-degree random regular graphs, macroscopic fluctuations are generally non-Gaussian and depend on the degree $d$~\cite[Theorem 19]{Johnson2015}.
By contrast, Theorem~\ref{thm:main} establishes universal Gaussian fluctuations, with limiting variance $V_g$ as in \eqref{eq: 77} (independent of $d$), within the restricted mesoscopic regime and test-function class considered here.
\end{remark}

\begin{remark}[Arbitrarily slowly decaying mesoscopic scales]
This case study shows that the universal mesoscopic CLT holds even when $\eta$ tends to zero arbitrarily slowly.
Mesoscopic CLTs for linear eigenvalue statistics often require a technical assumption on the decay rate of $\eta$, for example, $\eta\le N^{-\eps}$ for an arbitrarily small but fixed $\eps>0$ \cite{MR3678478,LS20,LLS24,Vova25a,EL26+}.
Notably, Peilen~\cite{MR4705298} established a mesoscopic CLT at arbitrary mesoscopic scales for $\beta$-ensembles.
\end{remark}

\begin{remark}[Comparison with correlated random matrix models]
In \cite{EL26+}, mesoscopic eigenvalue fluctuations are shown to be universal for a class of correlated random matrices.
A key assumption is the decay of correlations.
As a representative example, one can consider a random matrix $X$ with centered entries of unit variance satisfying
\begin{equation*}
	|\E[X_{ij}X_{kl}]|\lesssim \frac{1}{1+\mathrm{d}(ij,kl)^{s}}, \quad s>2,
\end{equation*}
where $\mathrm{d}(ij,kl)\coloneqq\min\{|i-k|+|j-l|,|i-l|+|j-k|\}$.
See \cite[Example 1.5]{EL26+} for more details.
Although this model allows long-range dependence, the decay of correlations suggests, at least heuristically, that its spectral behavior should resemble that of Wigner matrices.
However, random regular graphs have a completely different dependence structure due to global degree constraints, which makes their spectral universality more intriguing.
\end{remark}

To prove Theorem~\ref{thm:main}, we adapt the approach of Johnson~\cite{Johnson2015} and Dumitriu--Johnson--Pal--Paquette~\cite{MR3078290} to our mesoscopic setting.
Using a Chebyshev series expansion, we reduce the eigenvalue problem to counting certain non-backtracking walks on the graph; see Section~\ref{sec: 545} for the precise formulation.
The leading contribution to these counts comes from walks that traverse cycles, possibly repeatedly.
Then, a Poisson approximation for the cycle counts reduces the problem to a weighted sum of independent Poisson random variables.

In the mesoscopic setting, we show that only single traversals of cycles contribute to the limiting fluctuations.
An analogous reduction appears in Johnson's growing-degree result~\cite[Theorem 21]{Johnson2015}.
Here, however, the degree remains fixed, and we justify this reduction through bounds on the Chebyshev coefficients of the rescaled test function; see Lemma~\ref{lem:coefficient-bounds} and Section~\ref{sec: 700}.
The main new ingredient is a precise asymptotic analysis of these coefficients and the identification of the universal limiting variance.
By expressing their leading behavior in terms of the Fourier transform of $g$, we obtain the variance formula in \eqref{eq: 77}; see Lemma~\ref{lem: 183} and Section~\ref{sec: 739}.

The restriction $\eta \gg (\log{N})^{-1}$ in Theorem~\ref{thm:main} arises from the reduction to cycle counts and the associated Poisson approximation; see \cite[Proposition 12 and Corollary 14]{Johnson2015}.
In our mesoscopic setting, the argument requires $r$ many Poisson approximations of walk counts where $r$ satisfies $r\eta\to\infty$.
However, the approximation error is bounded by $O(\sqrt{r}(d-1)^{3r/2}/N)$.
When $\eta=N^{-\eps}$ for a small fixed $\eps>0$, the condition $r\eta\to\infty$ makes this error bound
blow up.
To establish the CLT throughout the full mesoscopic regime, one possible approach is to use resolvent methods building on recent advances in local laws and microscopic eigenvalue universality \cite{HuangMcKenzieYau2024,BH26+,HeHuangWang2026}.

\begin{remark}[Relation to earlier methods]
Earlier graph-theoretic aspects of the connection between orthogonal polynomials and non-backtracking walks were studied by McKay~\cite{McKay1983} and by Li and Sol\'e~\cite{LiSole1996}.
This connection also underlies the modified moment method developed systematically by Sodin~\cite{Sodin2007}.
This approach was applied to spectral-edge universality by Feldheim and Sodin~\cite{FeldheimSodin2010} and to the spectral
edges of random band matrices by Sodin~\cite{Sodin2010}.
Our variance calculation is related to earlier Fourier-based analyses of mesoscopic fluctuations \cite{ErdosKnowles2015I,ErdosKnowles2015II,MR1797877,MR4216520}.
In our result, the universal variance emerges from the mesoscopic asymptotics of the Chebyshev coefficients arising from the
cycle-count representation.
\end{remark}

\begin{remark}[Other related results]
For random regular graphs with sufficiently slowly growing degree, Gaussian fluctuations of macroscopic linear eigenvalue statistics were established in \cite{MR3078290,Johnson2015}.
Johnson and Pal~\cite{JohnsonPal2014} studied joint fluctuations of polynomial linear eigenvalue statistics for sequentially growing random regular multigraphs in the permutation model.
Bulk and edge universality have also been proved in various growing-degree regimes \cite{BauerschmidtHuangKnowlesYau2017,HuangYau2026Growing,BH26+}.
Recent progress on eigenvectors includes Gaussian-wave limits for fixed-degree random regular graphs at the spectral edge
\cite{HeHuangYau2025} and in the bulk \cite{HeHuangWang2026}.

For sparse Erd\H{o}s--R\'enyi graphs, global CLTs for linear eigenvalue statistics with constant expected degree were proved
by Shcherbina and Tirozzi~\cite{ShcherbinaTirozzi2010}.
In regimes with polynomially growing expected degree, He~\cite{He2020Sparse} established mesoscopic CLTs for resolvent
traces of centered sparse matrices, exhibiting a transition between universal GOE fluctuations and sparsity-dependent Gaussian fluctuations.
Recent results also establish asymptotic joint Gaussianity of nontrivial edge eigenvectors \cite{HeHuangWang2026Eigenvectors}.

Beyond these models, global spectral fluctuation results have been obtained for inhomogeneous random graphs with graphon limits \cite{ZhuZhu2024}.
Hirsch--Nam--Otto~\cite{HirschNamOtto2026} proved CLTs for linear eigenvalue statistics of random geometric graphs, with extensions to $k$-nearest neighbor graphs and relative neighborhood graphs for polynomial test functions.
For two ensembles of random quantum graphs, including one whose underlying graph is a random regular graph, Maltsev and Osman~\cite{MO26+} established mesoscopic variance asymptotics matching the GOE/GUE predictions.
\end{remark}

\subsection*{Organization}
Section \ref{sec: 209} describes the proof of Theorem \ref{thm:main} assuming several technical lemmas.
These lemmas are proved in Section \ref{sec: 331}.

\subsection*{Notation}
All limits are taken as $N\to\infty$ with $d$, $E$, and $g$ fixed.
We occasionally suppress the dependence on $N$ for brevity.
The letters $C,c>0$ denote constants that may change from line to line.
Unless otherwise stated, these constants may depend on $d,E,g$ and fixed auxiliary parameters, but are independent of $N$ and the summation indices.
For nonnegative quantities $a,b$, we write $a\lesssim b$ if
$a\le Cb$, and $a\sim b$ if both $a\lesssim b$ and $b\lesssim a$.
For positive sequences $a_N,b_N$, the notations $a_N\ll b_N$,
$b_N\gg a_N$, and $a_N=o(b_N)$ mean that $a_N/b_N\to0$.
We write $a_N=O(b_N)$ if $|a_N|\lesssim b_N$.

We use $\Longrightarrow$ for convergence in distribution and $\mathcal{N}(m,v)$ for the Gaussian distribution with mean $m$ and variance $v$.
For random elements $X,Y$ taking values in the same measurable space,
we denote the total variation distance between their laws by
\begin{equation*}
    \TV(X,Y)
    \coloneqq\sup_B\bigl|\mathbb{P}(X\in B)-\mathbb{P}(Y\in B)\bigr|,
\end{equation*}
where the supremum is taken over all measurable sets.

\subsection*{Acknowledgements}
This work is supported by ERC Advanced Grant ``RMTBeyond'' No.\,101020331.
The author is grateful to L\'{a}szl\'{o}~Erd\H{o}s for helpful discussions.

\section{Proof of Theorem \ref{thm:main}}\label{sec: 209}

Define $f\equiv f_N$ by
\begin{equation*}
	f(x) \coloneqq g\left(\frac{x-E}{\eta}\right),
\end{equation*}
which is the rescaled version of $g$.
Let $\{T_{k}\}_{k\ge 0}$ be the \emph{Chebyshev polynomials} of the first kind:
\begin{equation}\label{eq: 111}
	T_{0}(x) = 1, \quad T_{1}(x) = x, \quad T_{n+1}(x) = 2x T_{n}(x) - T_{n-1}(x) \;\; (n\ge 1).
\end{equation}
The following lemma provides the Chebyshev series of $f$.

\begin{lemma}[Chebyshev series]\label{lem: 108}
We have
\begin{equation}\label{eq:cheb-expansion}
     f(x) = \sum_{k=0}^{\infty}c_{k,N}T_k(x/2),
\end{equation}
where
\begin{equation}\label{eq:coefficients}
     c_{k,N} = \frac{2}{\pi}\int_{0}^{\pi} g\left(\frac{2\cos{\theta}-E}{\eta}\right) \cos(k\theta)\,\mathrm{d}\theta,
     \quad k\geq 1,
\end{equation}
and $c_{0,N} = \frac{1}{\pi}\int_{0}^{\pi}g\big(\frac{2\cos{\theta}-E}{\eta}\big)\,\mathrm{d}\theta$.
The series \eqref{eq:cheb-expansion} converges absolutely and locally uniformly.
\end{lemma}

Now we set the truncation cutoff
\begin{equation}\label{eq:cutoff}
   r_N \coloneqq \left\lfloor\frac{\log N}{3\log(d-1)}\right\rfloor,
\end{equation}
and define the truncation
\begin{equation*}
	p(x) \coloneqq \sum_{k=0}^{r_N}c_{k,N}T_k(x/2).
\end{equation*}

\begin{lemma}[Truncation]\label{lem:trace-truncation}
The following holds deterministically:
\begin{equation*}
	\left| \bigl(\Tr f(H)-\E\Tr f(H)\bigr) - \bigl(\Tr p(H)-\E\Tr p(H)\bigr)\right|
	\longrightarrow 0.
\end{equation*}
\end{lemma}

By Lemma \ref{lem:trace-truncation}, it suffices to show that
\begin{equation}\label{eq: 139}
	\Tr p(H)-\E\Tr p(H) \Longrightarrow\ \mathcal{N}(0,V_{g}).
\end{equation}
We introduce independent Poisson random variables $Z_j$ ($j\geq3$) with means
\begin{equation}\label{eq: 445}
	\lambda_j:=\E Z_j=\frac{(d-1)^j}{2j},
\end{equation}
and define
\begin{equation}\label{eq:poisson-walks}
    Q_k := \sum_{\substack{j\mid k\\j\geq3}}2jZ_j.
\end{equation}

\begin{lemma}[Poisson approximation]\label{lem: 153}
In order to show \eqref{eq: 139}, it is enough to prove
\begin{equation}\label{eq: 153}
	\sum_{k=3}^{r_N}\frac{c_{k,N}}{2(d-1)^{k/2}}\left(Q_k-\E Q_k\right) \Longrightarrow\ \mathcal{N}(0,V_{g}).
\end{equation}
\end{lemma}

Using \eqref{eq:poisson-walks} and reordering,
\begin{equation}\label{eq: 175}
	\sum_{k=3}^{r_N}\frac{c_{k,N}}{2(d-1)^{k/2}}\left(Q_k-\E Q_k\right)
	= \sum_{j=3}^{r_N} j(Z_j-\lambda_j)\sum_{1\le m\le r_{N}/j}c_{mj,N}(d-1)^{-mj/2}.
\end{equation}
The following lemma shows that the term with $m=1$,
\begin{equation*}
	\sum_{j=3}^{r_N}j c_{j,N}(d-1)^{-j/2}(Z_j-\lambda_j),
\end{equation*}
has the leading contribution.

\begin{lemma}[Single traversals of cycles]\label{lem: 169}
It follows that
\begin{equation*}
	\lim_{N\to\infty}\E\left| \sum_{k=3}^{r_N}\frac{c_{k,N}}{2(d-1)^{k/2}}\left(Q_k-\E Q_k\right) - \sum_{j=3}^{r_N}j c_{j,N}(d-1)^{-j/2}(Z_j-\lambda_j) \right|^{2} = 0.
\end{equation*}
\end{lemma}

Applying Lemma \ref{lem: 169}, to show \eqref{eq: 153}, it is sufficient to check
\begin{equation}\label{eq: 177}
	Y_N \coloneqq \sum_{j=3}^{r_N}j c_{j,N}(d-1)^{-j/2}(Z_j-\lambda_j) \Longrightarrow\ \mathcal{N}(0,V_{g}).
\end{equation}
We indeed establish \eqref{eq: 177} by the following lemma.

\begin{lemma}[Gaussianity and variance computation]\label{lem: 183}
For each fixed $t\in\R$, we have
\begin{equation}\label{eq: 203}
	\log\E e^{itY_N} = -\frac{t^2S_N}{2}+O(t^{3}\eta^3),
\end{equation}
where
\begin{equation*}
	S_N \coloneqq \Var Y_N 
	=\frac12\sum_{j=3}^{r_N}j c_{j,N}^2.
\end{equation*}
Furthermore,
\begin{equation}\label{eq: 213}
	\lim_{N\to\infty} S_N = V_g,
\end{equation}
where $V_{g}$ is defined in \eqref{eq: 77}.
\end{lemma}

Lemma \ref{lem: 183} implies that the limit of $Y_N$ is Gaussian with the desired limiting variance.
Therefore, combining Lemmas \ref{lem: 108}--\ref{lem: 183}, we complete the proof of Theorem \ref{thm:main}.
\qed

\bigskip
\noindent
We will prove Lemmas \ref{lem: 108}--\ref{lem: 183} in the next section.

\section{Proof of lemmas}\label{sec: 331}

\subsection{Proof of Lemmas \ref{lem: 108}--\ref{lem:trace-truncation}: Chebyshev series and truncation}

\begin{proof}[Proof of Lemma \ref{lem: 108}]

Note that $T_k(\cos\theta)=\cos(k\theta)$.
For each fixed $N$, the Chebyshev series and coefficient formula on $[-2,2]$ follow from \cite[Theorem~3.1]{Trefethen2013}.
Thus it is enough to justify the Chebyshev series outside $[-2,2]$ for completeness.

Fix $N$. Recall that $\widehat{g}$ is compactly supported by assumption.
Then, we can choose $B>0$ such that $\operatorname{supp}\widehat g\subset[-B,B]$.
Since $g$ is Schwartz, Fourier inversion gives the extension of $g$ to the complex plane:
\begin{equation}\label{eq: 240}
    g(z)=\frac{1}{2\pi}\int_{-B}^{B}
        \widehat g(t)e^{itz}\,\mathrm{d}t,
    \qquad z\in\mathbb C,
\end{equation}
which defines an entire function and agrees with the original $g$ on the real line.
Next define
\begin{equation}\label{eq: 248}
    F_N(w) \coloneqq g\left(\frac{w+w^{-1}-E}{\eta}\right),
\end{equation}
which is holomorphic on $\mathbb C\setminus\{0\}$. Its Laurent series
\begin{equation}\label{eq: 252}
    F_N(w)=\sum_{m\in\mathbb Z}a_mw^m
\end{equation}
converges absolutely and locally uniformly throughout $\mathbb C\setminus\{0\}$.
We calculate each Laurent coefficient. Set
\begin{equation*}
    h_N(\theta) \coloneqq g\left(\frac{2\cos\theta-E}{\eta}\right).
\end{equation*}
For $w=e^{i\theta}$, we have $w+w^{-1}=2\cos\theta$, so $F_N(e^{i\theta})=h_N(\theta)$. 
Applying the Laurent coefficient formula on the unit circle, for every $m\in\mathbb Z$,
\begin{equation*}
	a_m = \frac{1}{2\pi i}\oint_{|w|=1} F_N(w)w^{-m-1}\,\mathrm{d}w
	= \frac{1}{2\pi}\int_{-\pi}^{\pi} h_N(\theta)e^{-im\theta}\,\mathrm{d}\theta.
\end{equation*}
Since $h_N$ is even, the sine terms vanish: for every $k\ge1$,
\begin{equation}\label{eq: 267}
    a_k = a_{-k} =\frac{1}{\pi}\int_0^\pi h_N(\theta)\cos(k\theta)\,\mathrm{d}\theta
    = \frac{c_{k,N}}{2},
\end{equation}
and $a_0 = \frac{1}{\pi}\int_0^\pi h_N(\theta)\,\mathrm{d}\theta = c_{0,N}$.
Hence, for every $w\ne0$,
\begin{equation*}
    F_N(w)
    =c_{0,N}
     +\sum_{k=1}^{\infty}
        \frac{c_{k,N}}{2}(w^k+w^{-k}).
\end{equation*}
We use the identity,
\begin{equation}\label{eq: 278}
    T_k\left(\frac{w+w^{-1}}{2}\right) = \frac{w^k+w^{-k}}{2}, \quad w\neq 0,
\end{equation}
to connect this Laurent series to Chebyshev polynomials.
Considering the functions $P_{k}(u)=(w^k+w^{-k})/2$ with $u=(w+w^{-1})/2$, we observe that $\{P_{k}\}_{k\ge 0}$ satisfies \eqref{eq: 111} and thus the identity \eqref{eq: 278} is obtained.

Given any $x\in\mathbb{R}$, choose a root $w$ of $w^2-xw+1=0$.
This root is nonzero and satisfies $w+w^{-1}=x$ for $w\neq 0$. Substituting it into the Laurent series proves
\begin{equation*}
	f(x) = g\left(\frac{x-E}{\eta}\right) = F_N(w) = c_{0,N} + \sum_{k=1}^{\infty}c_{k,N}T_k(x/2).
\end{equation*}
Thus the Chebyshev series holds for every $x\in\mathbb{R}$, in particular outside $[-2,2]$.
By the convergence of the Laurent series of $F_N$, its coupled Chebyshev series also converges absolutely and locally uniformly outside $[-2,2]$ because for $|x|>2$ we can choose $x=\operatorname{sgn}(x)\,2\cosh t$ with some $t>0$ and set $w=\operatorname{sgn}(x)e^t$ so that $w+w^{-1}=x$.

\end{proof}

To prove Lemma \ref{lem:trace-truncation} and also for later usage, we need the following estimates.

\begin{lemma}\label{lem:coefficient-bounds}
There is a constant $C$, independent of $k$ and $N$, such that
\begin{equation}\label{eq:uniform-coeff}
    |c_{k,N}|\leq C\eta. 
\end{equation}
Choose $B>0$ such that $\supp\widehat{g}\subset[-B,B]$. Then, whenever $k\eta>B$,
\begin{equation}\label{eq:fast-coeff}
    |c_{k,N}| \leq C\left(\frac{eB}{k\eta}\right)^k.
\end{equation}
\end{lemma}

\begin{proof}
Set $\theta_{E}=\arccos(E/2)$.
First we shall choose a sufficiently small (fixed) neighborhood $I(\theta_{E})$ of $\theta_{E}$ such that
$\sin\theta$ is bounded away from zero in $I(\theta_{E})$.
This is possible because $E\in(-2,2)$ in our setting.
Now we decompose the integral \eqref{eq:coefficients} into two parts: (i) integral over $I(\theta_{E})$ and (ii) integral over  $[0,\pi]\backslash I(\theta_{E})$.
For the integral over $I(\theta_{E})$, using the change of variables $t=(2\cos\theta-E)/\eta$,
\begin{equation*}
	\bigg|\int_{I(\theta_{E})}g\left(\frac{2\cos{\theta}-E}{\eta}\right) \cos(k\theta)\,\mathrm{d}\theta\bigg|
	\lesssim \eta,
\end{equation*}
where we used that $|\sin\theta|=|\frac{\mathrm{d}}{\mathrm{d}\theta}\cos\theta|$ is away from zero in $I(\theta_{E})$.
For the integral over $[0,\pi]\backslash I(\theta_{E})$, $|2\cos\theta-E|$ is bounded away from zero.
Then, by the Schwartz decay of $g$,
\begin{equation*}
	\bigg|\int_{[0,\pi]\backslash I(\theta_{E})}g\left(\frac{2\cos{\theta}-E}{\eta}\right) \cos(k\theta)\,\mathrm{d}\theta\bigg|\lesssim \eta.
\end{equation*}
This proves \eqref{eq:uniform-coeff}, uniformly in $k$ and $N$.

Next we shall prove the second bound \eqref{eq:fast-coeff}.
Recall that $\supp\widehat{g}\subset[-B,B]$ and the Fourier inversion \eqref{eq: 240} extends $g$ to an entire function.
In particular, since $|e^{itz}|\le e^{B|\operatorname{Im}z|}$ for $t\in[-B,B]$, we get
\begin{equation}\label{eq: 144}
    |g(z)| \lesssim e^{B|\operatorname{Im}z|}.
\end{equation}
Recall $F_N$ defined in \eqref{eq: 248} and its Laurent series \eqref{eq: 252}.
By \eqref{eq: 267}, we already find that $a_{k}=c_{k,N}/2$ for every $k\ge 1$.
Thus, in order to get a bound for $c_{k,N}$, we shall consider the Laurent coefficient $a_k$.
Choose $\rho>1$. 
The function $F_N(w)w^{-k-1}$ is holomorphic in the annulus $1\le |w|\le\rho$.
We can deform the unit circle to the circle of radius $\rho$:
\begin{equation}\label{eq: 178}
    c_{k,N}
    =\frac{1}{\pi i}\oint_{|w|=\rho}
        F_N(w)w^{-k-1}\,\mathrm{d}w, \quad k\ge 1.
\end{equation}
To bound $F_N$ on this circle, write $w=\rho e^{i\theta}$. Then
\begin{equation*}
    w+w^{-1}
    =(\rho+\rho^{-1})\cos\theta
      +i(\rho-\rho^{-1})\sin\theta,
\end{equation*}
so, because $E$ is real,
\begin{equation*}
    \left|
    \operatorname{Im} \bigg(\frac{w+w^{-1}-E}{\eta}\bigg)
    \right|
    \le \frac{\rho-\rho^{-1}}{\eta}.
\end{equation*}
Consequently, the bound \eqref{eq: 144} on $g$ implies
\begin{equation}\label{eq: 197}
    \max_{|w|=\rho}|F_N(w)|
    \lesssim \exp\left(\frac{B(\rho-\rho^{-1})}{\eta}\right).
\end{equation}
In conclusion, it follows from \eqref{eq: 178} and \eqref{eq: 197}, for every $\rho>1$,
\begin{equation*}
    |c_{k,N}|
    \le \frac{1}{\pi}(2\pi\rho)\rho^{-k-1} \max_{|w|=\rho}|F_N(w)|
    \lesssim \exp\left(\frac{B(\rho-\rho^{-1})}{\eta}\right)\rho^{-k}.
\end{equation*}
Finally, if $k\eta>B$, we set $\rho=k\eta/B>1$. Simply with $\rho-\rho^{-1}\le\rho$, we obtain
\begin{equation*}
    |c_{k,N}| \lesssim \exp\left(\frac{B\rho}{\eta}\right)\rho^{-k} = \left(\frac{eB}{k\eta}\right)^k,
\end{equation*}
which shows \eqref{eq:fast-coeff}.
\end{proof}

\begin{proof}[Proof of Lemma \ref{lem:trace-truncation}]

Our assumption $\eta\log N\to\infty$ implies (recalling that $d$ is fixed)
\begin{equation}\label{eq:r-eta}
r_N\eta\longrightarrow\infty.
\end{equation}
Since $A$ has nonnegative entries and every row sums to $d$, every eigenvalue of $A$ has absolute value at most $d$.
Every eigenvalue of $H_N$ lies in $[-d/\sqrt{d-1},d/\sqrt{d-1}]$ and thus it is enough to show
\begin{equation}\label{eq:uniform-truncation}
N\sup_{|x|\leq d/\sqrt{d-1}}|f(x)-p(x)|
\longrightarrow0.
\end{equation}
Note that
\begin{equation*}
	\sup_{|x|\leq d/\sqrt{d-1}}|T_k(x/2)|\leq (d-1)^{k/2},
\end{equation*}
where we used \eqref{eq: 278}
and $T_{k}(\pm\cosh{t})=(\pm 1)^{k}\cosh(kt)$ with $d/\sqrt{d-1}=\sqrt{d-1}+1/\sqrt{d-1}=2\cosh(a)$ for $a=\log(\sqrt{d-1})$.
Introduce
\begin{equation*}
	b_N \coloneqq \frac{eB\sqrt{d-1}}{(r_N+1)\eta}.
\end{equation*}
By \eqref{eq:r-eta}, we see that $b_N\to0$. For all $k>r_N$ and sufficiently large $N$, the estimate \eqref{eq:fast-coeff} in Lemma \ref{lem:coefficient-bounds} implies $|c_{k,N}|(d-1)^{k/2}\lesssim b_N^k$. As a result,
\begin{equation*}
	\sup_{|x|\leq d/\sqrt{d-1}}|f(x)-p(x)| \lesssim \sum_{k>r_N}b_N^k = \frac{b_N^{r_N+1}}{1-b_N}.
\end{equation*}
Now $r_N/\log N\to(3\log(d-1))^{-1}>0$, whereas $\log b_N\to-\infty$.
Therefore $N b_N^{r_N+1}\to0$ and we obtain \eqref{eq:uniform-truncation} and complete the proof of Lemma \ref{lem:trace-truncation}.

\end{proof}

\subsection{Proof of Lemma \ref{lem: 153}: Poisson approximation}\label{sec: 545}

\begin{figure}[h]
\begin{center}
\begin{tikzpicture}[
  vertex/.style={circle,draw,fill=white,inner sep=0pt,
                 minimum size=5pt,outer sep=1pt},
  graph/.style={gray!55,line width=0.8pt},
  walk/.style={-{Stealth[length=2.2mm]},line width=1.1pt},
  walklabel/.style={font=\small,fill=white,inner sep=1.5pt},
  every label/.style={font=\small}
]

\begin{scope}
  \node[vertex,label=left:$a$]  (a) at (0,0) {};
  \node[vertex,label=below:$b$] (b) at (1.6,0) {};
  \node[vertex,label=above:$c$] (c) at (3.1,1) {};
  \node[vertex,label=below:$d$] (d) at (3.1,-1) {};

  \draw[graph] (a)--(b)--(c)--(d)--(b);

  \draw[walk,red!75!black] (a) to[bend left=22]
    node[walklabel,above] {$1$} (b);
  \draw[walk,red!75!black] (b) to[bend left=22]
    node[walklabel,below] {$5$} (a);
  \draw[walk,red!75!black] (b)--node[walklabel,above left] {$2$} (c);
  \draw[walk,red!75!black] (c)--node[walklabel,right] {$3$} (d);
  \draw[walk,red!75!black] (d)--node[walklabel,below left] {$4$} (b);
\end{scope}

\begin{scope}[xshift=6cm]
  \node[vertex,label=left:$a$]  (a) at (0,0) {};
  \node[vertex,label=below:$b$] (b) at (1.6,0) {};
  \node[vertex,label=above:$c$] (c) at (3.1,1) {};
  \node[vertex,label=below:$d$] (d) at (3.1,-1) {};

  \draw[graph] (a)--(b)--(c)--(d)--(b);

  \draw[walk,blue!75!black] (b)--node[walklabel,above left] {$1$} (c);
  \draw[walk,blue!75!black] (c)--node[walklabel,right] {$2$} (d);
  \draw[walk,blue!75!black] (d)--node[walklabel,below left] {$3$} (b);
\end{scope}
\end{tikzpicture}
\end{center}
\caption{Two closed non-backtracking walks on the same underlying graph.
The left walk is not cyclically non-backtracking, since its last step
$b\to a$ reverses its first step $a\to b$; the right walk is cyclically
non-backtracking. Numbers indicate the order of traversal, and the two
curved arrows on the left represent traversals of the same edge $ab$.}\label{fig: 493}
\end{figure}
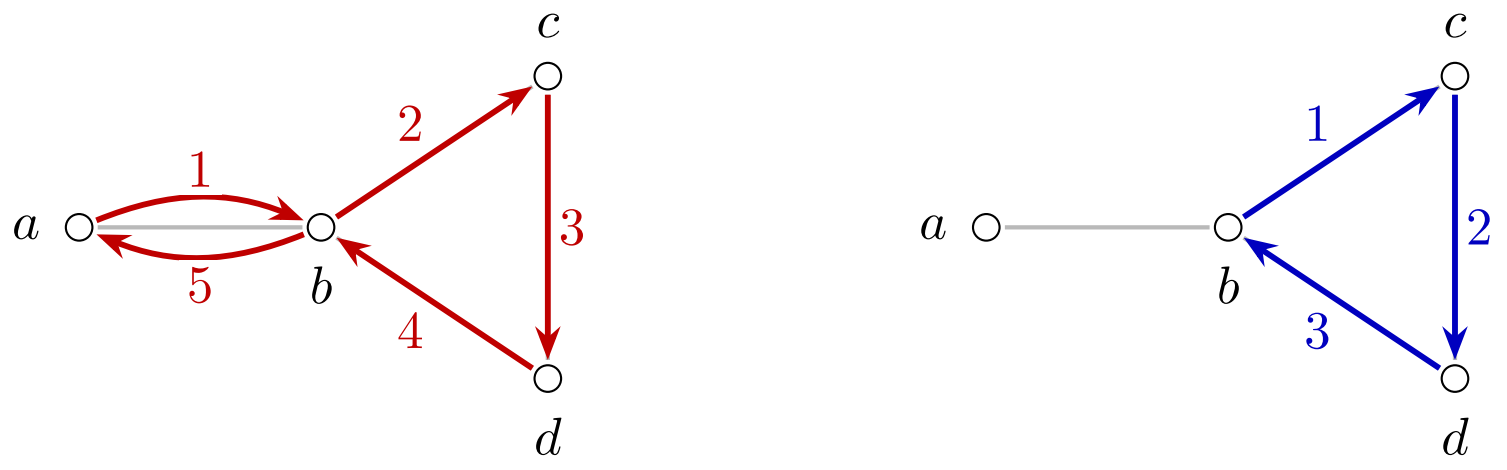

As in \cite{Johnson2015,MR3078290,MR2437174}, we introduce some notions of walks on a graph:
\begin{equation*}
	\text{\emph{non-backtracking walk} and \emph{cyclically non-backtracking walk}.}
\end{equation*}
A non-backtracking walk is a walk that never traverses an edge and then immediately traverses the same edge in the reverse direction.
A closed non-backtracking walk is called a cyclically non-backtracking walk if its last step is not the reverse of its first step.
See Figure \ref{fig: 493} below for example.

Let $W_k^{(N)}$ denote the number of cyclically non-backtracking walks of length $k$ in $\mathcal{G}_N$.
In \cite{Johnson2015,MR3078290}, the following polynomials are introduced:
\begin{equation*}
	\Gamma_{0}(x)=1, \quad \Gamma_{2m}(x)=2T_{2m}(x/2)+\frac{d-2}{(d-1)^m}, \quad 
	\Gamma_{2m-1}(x)=2T_{2m-1}(x/2)\;\;(m\geq1). 
\end{equation*}
According to \cite[Proposition 15]{Johnson2015} (also see \cite[Proposition 32]{MR3078290}), 
we have
\begin{equation}\label{eq:trace-walk}
\Tr\Gamma_k(H)=(d-1)^{-k/2}W_k^{(N)}.
\end{equation}
The difference between $\Gamma_k$ and $2T_k(x/2)$ is a deterministic constant, which is cancelled out by centering.
Note that $W_1^{(N)}=W_2^{(N)}=0$ in a simple graph.
We now have
\begin{equation}\label{eq:centered-polynomial}
\Tr p_N(H_N)-\E\Tr p_N(H_N)
=
\sum_{k=3}^{r_N}
\frac{c_{k,N}}{2(d-1)^{k/2}}
\left(W_k^{(N)}-\E W_k^{(N)}\right).
\end{equation}

We also introduce independent Poisson random variables $Z_j$ ($j\geq3$) with means
\begin{equation}
	\lambda_j:=\E Z_j=\frac{(d-1)^j}{2j},
\end{equation}
and define
\begin{equation}
    Q_k := \sum_{\substack{j\mid k\\j\geq3}}2jZ_j.
\end{equation}
A simple cycle of length $j$ contributes $2j$ cyclically non-backtracking walks of length $k$ when $j\mid k$ (by the choice of starting point and direction).
The variables $Q_k$ are thus the walk counts constructed from independent Poisson cycle counts.
It follows from \cite[Corollary 14]{Johnson2015} that
\begin{equation}\label{eq:tv-bound}
    \TV\left((W_k^{(N)})_{3\leq k\leq r_N},(Q_k)_{3\leq k\leq r_N}\right)
    \lesssim \frac{\sqrt{r_N}\,(d-1)^{3r_N/2-1}}{N}.
\end{equation}
Since $(d-1)^{r_N}\leq N^{1/3}$, the resulting bound would be $O(\sqrt{\log N}\,N^{-1/2})=o(1)$.

Next we should check the required expectation comparison: $|\E W_k^{(N)}-\E Q_k|$. 
Let $C_j^{(N)}$ be the number of simple cycles of length $j$.
According to \cite[Proposition 12]{Johnson2015},
\begin{equation}\label{eq: 465}
	B_k^{(N)} \coloneqq W_k^{(N)} - \sum_{\substack{j\mid k\\j\geq3}}2jC_j^{(N)}, \quad
	\E B_k^{(N)}\lesssim \frac{k^6(d-1)^k}{N},
\end{equation}
where $B_k^{(N)}$ counts walks that are not repeated walks over a cycle. 
By McKay--Wormald--Wysocka \cite[Eq\,(2.2)]{MWW2004},
\begin{equation}\label{eq:cycle-mean}
    \E C_j^{(N)} = \frac{(d-1)^j}{2j}\left(1+O\left(\frac{j(j+d)}{N}\right)\right),\qquad 3\leq j\leq r_N.
\end{equation}
Combining \eqref{eq: 445}--\eqref{eq:poisson-walks} and \eqref{eq: 465}--\eqref{eq:cycle-mean},
\begin{equation}\label{eq:mean-discrepancy}
    \left|\E W_k^{(N)}-\E Q_k\right| \lesssim \frac{\big(k^6+k(k+d)\big)(d-1)^k}{N}.
\end{equation}
Multiply \eqref{eq:mean-discrepancy} by the coefficient in
\eqref{eq:centered-polynomial}, sum, and use
\eqref{eq:uniform-coeff}. The total difference between the two
centers is bounded by
\begin{align}
\sum_{k=3}^{r_N}
\frac{|c_{k,N}|}{2(d-1)^{k/2}}
\left|\E W_k^{(N)}-\E Q_k\right|
&\lesssim \frac{\eta}{N}
\sum_{k=3}^{r_N}\big(k^6+k(k+d)\big)(d-1)^{k/2}\nonumber\\
&\leq
C\eta \,r_N^7\frac{(d-1)^{r_N/2}}{N}\nonumber\\
&\leq
C\eta(\log N)^7 N^{-5/6}
\longrightarrow0.\label{eq: 495}
\end{align}
Applying the total variation estimate \eqref{eq:tv-bound}
together with \eqref{eq:centered-polynomial} and \eqref{eq: 495}, we complete the proof of Lemma \ref{lem: 153}. 
\qed

\subsection{Proof of Lemma \ref{lem: 169}: contributions from single traversals over cycles}\label{sec: 700}

Define
\begin{equation}\label{eq:poisson-statistic}
    P_N := \sum_{k=3}^{r_N}\frac{c_{k,N}}{2(d-1)^{k/2}}\left(Q_k-\E Q_k\right).
\end{equation}
Recalling \eqref{eq: 175},
\begin{equation*}
	P_N
	= \sum_{j=3}^{r_N} j (Z_j-\lambda_j) \sum_{1\le m\le r_{N}/j}c_{mj,N}(d-1)^{-mj/2}.
\end{equation*}
The term corresponding to $m=1$ represents the Poisson approximation for single traversals of cycles.
Its sum is denoted by $Y_N$ in \eqref{eq: 177}:
\begin{equation*}
    Y_N = \sum_{j=3}^{r_N}j c_{j,N}(d-1)^{-j/2}(Z_j-\lambda_j).
\end{equation*}
We denote the remaining terms ($m\geq2$) by
\begin{equation*}
	R_N:=P_N-Y_N.
\end{equation*}
Using the estimate \eqref{eq:uniform-coeff} in Lemma \ref{lem:coefficient-bounds}, 
for each $j\geq3$,
\begin{equation*}
	\left|\sum_{\substack{m\geq2\\mj\leq r_N}}c_{mj,N}(d-1)^{-mj/2}\right|\lesssim \frac{\eta(d-1)^{-j}}{1-(d-1)^{-j/2}}
	\lesssim \eta (d-1)^{-j}.
\end{equation*}
Since the Poisson random variables $Z_{j}$ are independent and $\Var Z_j=\lambda_j=(d-1)^j/(2j)$,
\begin{equation*}
	\Var R_N \lesssim \eta^2\sum_{j=3}^{\infty} j (d-1)^{-j} = O(\eta^2) \longrightarrow0,
\end{equation*}
which means $P_N-Y_N\to0$ in $L^2$.
Therefore Lemma \ref{lem: 169} holds.
\qed

\subsection{Proof of Lemma \ref{lem: 183}: Gaussianity and variance computation}\label{sec: 739} 

The variance of the cycle sum is
\begin{equation}\label{eq:prelimit-variance}
    S_N = \Var Y_N = \sum_{j=3}^{r_N}j^2c_{j,N}^2 (d-1)^{-j}\frac{(d-1)^j}{2j}=\frac12\sum_{j=3}^{r_N}j c_{j,N}^2.
\end{equation}
Plus, the third order cumulant of $Y_{N}$ is bounded by
\begin{equation*}
\sum_{j=3}^{r_N}\lambda_j\left|j c_{j,N}(d-1)^{-j/2}\right|^3
=\frac{1}{2}\sum_{j=3}^{r_N}j^2|c_{j,N}|^3 (d-1)^{-j/2}
\lesssim \eta^3\sum_{j=3}^{\infty}j^2 (d-1)^{-j/2}
= O(\eta^3).
\end{equation*}
Since every cumulant of the Poisson distribution is equal to its mean, higher order cumulants are bounded analogously.
Thus, for each fixed $t\in\R$, considering the cumulant generating function, we obtain \eqref{eq: 203}.

\bigskip

Write
\begin{equation*}
	\theta_E:=\arccos(E/2), \quad v:=\sqrt{4-E^2}>0.
\end{equation*}
As in Lemma \ref{lem:coefficient-bounds}, recall that we choose $B>0$ satisfying $\supp\widehat{g}\subset[-B,B]$.

\begin{lemma}\label{lem:local-coeff}
For each fixed $M>1$, uniformly over integers $k\geq1$ with
$k\eta\leq M$,
\begin{equation}\label{eq:local-coeff}
c_{k,N}
=
\frac{2\eta}{\pi v}
\operatorname{Re}\!\left[
e^{ik\theta_E}
\widehat g\!\left(\frac{k\eta}{v}\right)
\right]
+O(M\eta^2).
\end{equation}
If $M>eB$, then
\begin{equation}\label{eq:variance-tail}
\sum_{k>M/\eta}k c_{k,N}^2\longrightarrow0.
\end{equation}
\end{lemma}

\begin{proof}
Let $\theta(x)=\arccos(x/2)$. Changing variables
$x=2\cos\theta$ in \eqref{eq:coefficients} gives
\begin{equation*}
	c_{k,N} =
	\frac2\pi\int_{-2}^{2}
	g\left(\frac{x-E}{\eta}\right)\frac{\cos(k\theta(x))}{\sqrt{4-x^2}}\,dx.
\end{equation*}
The part with $x$ outside a small (fixed) neighborhood of $E$ is negligible to any power of $\eta$
by the Schwartz decay of $g$.
Inside that neighborhood, write $x=E+\eta u$. By Taylor's theorem,
\begin{equation*}
	\theta(E+\eta u) = \theta_E-\frac{\eta u}{v} +O(\eta^2u^2), \qquad
    \frac1{\sqrt{4-(E+\eta u)^2}} = \frac1v+O(\eta|u|).
\end{equation*}
When $k\eta\leq M$, the error in $\cos(k\theta(x))$ is $O(M\eta u^2)$.
Since $g$ has finite moments of every order (e.g.~$\int_{\R}u^2|g(u)|\,\mathrm{d}u<\infty$),
by integrating with $|g(u)|$ (recalling $dx = \eta du$), we have a uniform $O(M\eta^2)$ error.
Then, the leading contribution of the integral would be
\begin{equation*}
	\frac{2\eta}{\pi v}\int_{\R}g(u)\cos\left(k\theta_E-\frac{k\eta}{v}u\right)\,du,
\end{equation*}
which equals the main term in \eqref{eq:local-coeff} by the definition of $\widehat g$.
Thus we obtain \eqref{eq:local-coeff}.

For $k\eta>M>eB$, the estimate \eqref{eq:fast-coeff} in Lemma \ref{lem:coefficient-bounds} implies
\begin{equation*}
	|c_{k,N}|\leq C(eB/M)^k.
\end{equation*}
Since $eB/M<1$ and $M/\eta\to\infty$, summing
$k(eB/M)^{2k}$ proves \eqref{eq:variance-tail}.
\end{proof}

Choose a fixed $M>eB$ large enough that $M/v>B + (\eta/v)$.
By $r_N\eta\to\infty$, the range $k\leq M/\eta$ lies below $r_N$ for all sufficiently large $N$.
Lemma \ref{lem:local-coeff} and \eqref{eq:prelimit-variance} show that
\begin{equation}\label{eq:variance-riemann}
   S_N = \frac{2}{\pi^2} \sum_{3 \le k \le M/\eta}
\frac{\eta}{v}\,t_k
\left[
\operatorname{Re}
\bigl(e^{ik\theta_E}\widehat g(t_k)\bigr)
\right]^2
+o(1), \quad t_k \coloneqq k\eta/v.
\end{equation}
Justification for the cumulative error is as follows.
Each error in \eqref{eq:local-coeff} is $O(M\eta^2)$ while its main term is $O(\eta)$.
Summing $k\,O(M\eta^3)$ over $k\leq M/\eta$, we have $O(M^{3}\eta)=o(1)$.

Next we use the following identity
\begin{equation*}
	\left[\operatorname{Re}(e^{ik\theta_E}z)\right]^2
	= \frac{|z|^2}{2} + \frac{1}{2}\operatorname{Re}(e^{2ik\theta_E}z^2).
\end{equation*}
The nonoscillatory part of \eqref{eq:variance-riemann},
\begin{equation*}
	\frac{1}{\pi^2} \sum_{3 \le k \le M/\eta}
    \frac{\eta}{v}\,t_k |\widehat{g}(t_k)|^2,
\end{equation*}
is a Riemann sum converging to
\begin{equation}\label{eq: 789}
	\frac1{\pi^2}\int_0^\infty t|\widehat g(t)|^2\,dt,
\end{equation}
where we used that our setting $M/v>B + (\eta/v)$ implies $t_{k} > B$ for $k=\lfloor M/\eta\rfloor$, and that $\widehat{g}(t)=0$ for $t>B$.
Thus the integral also stops at $B$.

We claim that the oscillatory part tends to zero:
\begin{equation}\label{eq: 795}
	\frac{1}{\pi^2} \sum_{3 \le k \le M/\eta} \frac{\eta}{v}\,t_k \operatorname{Re}(e^{2ik\theta_E}\widehat{g}(t_k)^2) = o(1).
\end{equation}
Indeed, $0<\theta_E<\pi$, so the partial sums of $e^{2ik\theta_E}$ are uniformly bounded. 
Summation by parts, applied to the smooth compactly supported function $t\widehat g(t)^2$, makes the corresponding Riemann sum $o(1)$. We prove the claim \eqref{eq: 795} in detail below for completeness.

\bigskip

Recall $E=2\cos\theta_E$ with $0<\theta_E<\pi$ and $v=\sqrt{4-E^2}$. Introduce
\begin{equation*}
    h_N \coloneqq \frac{\eta}{v},
    \quad
    F(t) \coloneqq t\,\widehat g(t)^2
    \quad (t\ge0).
\end{equation*}
Since $E$ is fixed and $\eta\to0$, we have $h_N\to0$.
Up to a fixed prefactor and taking a real part, the oscillatory
contribution is of the form
\begin{equation*}
    \widetilde{S}_N \coloneqq h_N \sum_{3 \le k \le M/\eta} e^{2ik\theta_E}F(kh_N).
\end{equation*}
It is enough to show that $\widetilde{S}_N\to0$. Set $z=e^{2i\theta_E}$ and
\begin{equation*}
    A_m=\sum_{k=3}^m z^k,
    \qquad A_2=0.
\end{equation*}
Because $0<\theta_E<\pi$, we have $z\ne1$. The geometric-series formula therefore gives
\begin{equation*}
    A_m=\frac{z^{3}(1-z^{m-2})}{1-z},
    \qquad
    |A_m|
    \le \frac{2}{|1-z|}
    =\frac{1}{\sin\theta_E}
    \quad\text{for every }m\ge3.
\end{equation*}
In particular, the bound on $A_m$ does not depend on $m$, and $\sin\theta_E$ is bounded away from zero.
Using $z^k=A_k-A_{k-1}$ and rearranging the finite sum, we obtain
the summation-by-parts identity
\begin{align*}
    \sum_{k=3}^{m}z^kF(kh_N)
    &=\sum_{k=3}^{m}
        (A_k-A_{k-1})F(kh_N)\\
    &=A_{m}F(mh_N)
      +\sum_{k=3}^{m-1}
        A_k\bigl(F(kh_N)-F((k+1)h_N)\bigr).
\end{align*}
Thus the bounded partial sums $A_k$ are multiplied by differences
between consecutive samples of $F$.
The assumption $\widehat g\in C_c^\infty(\mathbb R)$ implies that
$F$ is smooth and compactly supported on $[0,\infty)$. In particular,
\begin{equation*}
    F'(t)
    =\widehat g(t)^2
      +2t\,\widehat g(t)\widehat g'(t),
    \qquad
    \int_0^\infty |F'(t)|\,\mathrm{d}t<\infty.
\end{equation*}
For every $k$, the fundamental theorem of calculus gives
\begin{equation*}
    |F((k+1)h_N)-F(kh_N)|
    \le
    \int_{kh_N}^{(k+1)h_N}
       |F'(t)|\,\mathrm{d}t.
\end{equation*}
The intervals of integration are disjoint. Combining this fact
with the summation-by-parts identity and the bound on $A_k$, we find
\begin{align}
    |\widetilde{S}_N|
    &\le
      \frac{h_N}{\sin\theta_E}
      \left(
        |F(\lfloor M/\eta \rfloor h_N)|
        +\sum_{3 \le k \le M/\eta}
          |F((k+1)h_N)-F(kh_N)|
      \right) \nonumber\\
    &\le
      \frac{h_N}{\sin\theta_E}
      \left(
        \|F\|_\infty
        +\int_0^\infty |F'(t)|\,\mathrm{d}t
      \right)
      =o(1). \label{eq: 878}
\end{align}
In fact, since $\lfloor M/\eta \rfloor h_N$ lies beyond the support of $F$, the boundary term $F(\lfloor M/\eta \rfloor h_N)$ is zero; this is not needed for the estimate though.
The estimate \eqref{eq: 878} applies to the complex sum $\widetilde{S}_N$, and hence also to its real part.
The claim \eqref{eq: 795} is proven.

\bigskip

From \eqref{eq:variance-riemann}, using the Riemann sum convergence to \eqref{eq: 789} and the claim \eqref{eq: 795},
we get
\begin{equation}\label{eq:fourier-variance}
S_N\longrightarrow
\frac1{\pi^2}\int_0^\infty
t|\widehat g(t)|^2\,dt.
\end{equation}
Since  $g$ is real-valued, it follows that $|\widehat g(-t)|=|\widehat g(t)|$.
Thus,
\begin{equation*}
	\int_0^\infty t|\widehat g(t)|^2\,dt = \frac{1}{2} \int_{\R}|t|\,|\widehat g(t)|^2\,dt.
\end{equation*}
Therefore we find that \eqref{eq:fourier-variance} implies \eqref{eq: 213}.
We complete the proof of Lemma \ref{lem: 183}. 
\qed

\section*{AI usage disclosure}
The author used an AI tool (ChatGPT; OpenAI; model: GPT-6.1 Sol Extra High) to assist with literature searches, language editing (including stylistic revisions and checks for typos and grammatical errors), preparation of the figure, and supplementary checks of proof details.
The author developed the main ideas and carried out the research.
All AI-assisted material was reviewed and verified by the author who takes full responsibility for the content of this paper.

\bibliographystyle{abbrv}

\end{document}